\documentclass[3p]{elsarticle}
\usepackage{amscd, amssymb, amsthm}
\usepackage[all]{xy}
\usepackage{color}
\usepackage{url}
\usepackage{amscd, amssymb, amsfonts, amsthm, amsmath}
\usepackage{mathdots}
\usepackage[breaklinks=true]{hyperref}
\usepackage{multicol}
\usepackage{yfonts}

\usepackage[normalem]{ulem}

\usepackage[labelsep=endash]{caption}

\numberwithin{equation}{section}
\usepackage[labelsep=endash]{caption}

\makeatletter
\def\ps@pprintTitle{%
  \let\@oddhead\@empty
  \let\@evenhead\@empty
  \let\@oddfoot\@empty
  \let\@evenfoot\@oddfoot
}
\makeatother

\newcommand{\topp}{\mbox{\rm{Top}}}
\newcommand{\length}{\mbox{\rm{l}}}

\newcommand{\start}{\mbox{\rm{s}}}
\newcommand{\target}{\mbox{\rm{t}}}
\newcommand{\gl}{\mbox{\rm{gldim}}}
\newcommand{\fin}{\mbox{\rm{findim}}}
\newcommand{\Fin}{\mbox{\rm{Findim}}}

\newcommand{\fidim}{\mbox{\rm{$\phi$dim}}}
\newcommand{\psidim}{\mbox{\rm{$\psi$dim}}}
\newcommand{\findim}{\mbox{\rm{findim}}}

\newcommand{\Dell}{\mbox{\rm{dell}}}
\newcommand{\Dellm}{\mbox{\rm{Dell}}}

\def\mod{\mbox{\rm{mod}}}
\def\Mod{\mbox{\rm{Mod}}}
\def\ind{\mbox{\rm{ind}}}
\def\add{\mbox{\rm{add}}}

\def\pd{\mbox{\rm{pd}}}
\def\id{\mbox{\rm{id}}}

\def\rk{\hbox{\rm{rk}}}

\def\Rep{\mbox{\rm{Rep}}}
\def\rep{\mbox{\rm{rep}}}

\def\tg{\mbox{\rm{t}}}
\def\supp{\mbox{\rm{supp}}}

\begin{document}
\newcommand{\mono}[1]{%
\gdef\puA{#1}}
\newcommand{\puA}{}
\newcommand{\faculty}[1]{%
\gdef\puC{#1}}
\newcommand{\puC}{}
\newcommand{\facultad}[1]{%
\gdef\puD{#1}}
\newcommand{\puD}{}
\newcommand{\N}{\mathbb{N}}
\newcommand{\Z}{\mathbb{Z}}
\newtheorem{teo}{Theorem}[section]
\newtheorem{prop}[teo]{Proposition}
\newtheorem{lema}[teo] {Lemma}
\newdefinition{ej}[teo]{Example}
\newtheorem{obs}[teo]{Remark}
\newtheorem{defi}[teo]{Definition}
\newtheorem{coro}[teo]{Corollary}
\newtheorem{nota}[teo]{Notation}



\title{Delooping levels}

\author{Marcos Barrios}
\ead{marcosb@fing.edu.uy}

\author{Marcelo Lanzilotta \corref{cor1}}
\ead{marclan@fing.edu.uy}

\author{Gustavo Mata \corref{cor2}}
\address{Universidad de la Rep\'ublica, Facultad de Ingenier\'ia -  Av. Julio Herrera y Reissig 565, Montevideo, Uruguay} 
\ead{gmata@fing.edu.uy}

\cortext[cor1]{The author was supported by the project of  MathAmSud - 23-MATH-06.}
\cortext[cor2]{Corresponding author}

\begin{abstract}{
} 
In \cite{G1} V. Gélinas introduced a homological invariant, called {\it delooping level} (\Dell), that bounds the finitistic dimension. In this article, we introduce another homological invariant (\Dellm) related to the delooping level for an Artin algebra. We compare this new tool with other dimensions as the finitistic dimension or the $\phi$-dimension (where $\phi$ is the first Igusa-Todorov function), and we also generalize Theorem 4.3. from \cite{G2} to truncated path algebras (Theorem \ref{TeoTruncada}). Finally, we show that for a monomial algebra $A$ the difference $\Dell(A) - \Fin(A)$ can be arbitrarily large (Example \ref{contraej}).
\end{abstract} 


\maketitle

\section{Introduction}
					
The finitistic dimensions are central invariants in the homological theory of Artin algebras. They are defined for $A$ an Artin algebra as follows:
\begin{itemize}
\item $\fin(A) = \{\pd (M) : M \in \mod A \text{ such that } \pd(M) < \infty\}$.
\item $\Fin(A) = \{\pd (M) : M \in \Mod A \text{ such that } \pd(M) < \infty\}.$
\end{itemize}

In 1960, H. Bass published two relevant questions about the finitistic dimensions (see \cite{Ba}) that have become conjectures over time:

\begin{itemize}
\item $\Fin(A) = \fin(A)$.
\item $\fin(A) < \infty$.
\end{itemize}

The first question is not true. In \cite{H-Z}, B. Huisgen-Zimmermann showed that there are monomial algebras such that $\Fin(A) = \fin(A)+1$. However, it was proved for radical square zero algebras (Theorem 3.1 from \cite{M}) and for truncated path algebras (Theorem 2.2 from \cite{BH-ZT}) that the conjecture is valid, and for monomial algebras it was proved that $\fin(A) \leq \Fin(A) \leq \fin(A)+1$ (Corollary B from \cite{H-Z}).\\

The delooping level for Noetherian semiperfect rings was introduced by V. Gélinas in \cite{G1}. In this article, he also proved the following relation between the finitistic dimension and the delooping level for every Noetherian semiperfect ring $A$:
\begin{equation}\label{eq1}
\Fin(A^{op}) \leq \Dell(A).
\end{equation}
In particular, the author proved that $\Fin(A^{op}) = \Dell(A)$ for every radical square zero algebra (Theorem 4.3. from \cite{G2}) and for every Gorenstein algebra (Example 2.17 from \cite{G1}). As a consequence of Theorem 3.1 from \cite{M}, we also have $\fin(A^{op}) = \Dell(A)$ for radical square zero algebras.

However, the inequality \ref{eq1} is not a proof for the finitistic dimension conjecture, since in \cite{KR} the authors showed that there exists an Artin algebra with infinite delooping level (still with finite finitistic dimension).\\

In an older attempt to prove the finitistic dimension conjecture, Igusa and Todorov introduced in \cite{IT} two functions which generalize the notion of projective dimension. Nowadays, they are known as the Igusa-Todorov functions, $\phi$ and $\psi$. The $\phi$-dimension and $\psi$-dimension (the supremum of $\phi$ and $\psi$ on $\mod A$ respectively) also verify for every Artin algebra $A$ that $\findim(A) \leq \fidim (A) \leq \psi \dim (A)$. These homological dimensions were studied in several works; we can cite, for example: \cite{BLMV}, \cite{BM1}, \cite{BM}, \cite{IT}, \cite{LM}, and \cite{LMM} among others.\\ 
  
In this article, we introduce another homological invariant related to the delooping level: Global delooping level (see Definition \ref{dell}, item 2). We study and compare this new invariant with other homological tools as the finitistic dimension or the $\phi$-dimension (Proposition \ref{sizigia finita}, Theorem \ref{n-Goreinstein}, and Proposition \ref{id finita}). We also generalize Theorem 4.3. from \cite{G2} to truncated path algebras (Theorem \ref{TeoTruncada}). Finally, we show that for a monomial algebra $A$, the difference $\Dell(A) - \Fin(A)$ can be arbitrary large (Example \ref{contraej}).

\section{Preliminaries}

Throughout this article, $A$ is an Artin algebra, $\Mod A$ is the category of right $A$-modules, and $\mod A$ is the full subcategory of $\Mod A$ determined the finitely generated right $A$-modules. $\ind A$ is the full subcategory of $\mod A$ formed by all the indecomposable modules, and $\mathcal{P}_A \subset \mod A$ is the class of projective $A$-modules. We use $\mathcal{S}_A$ for the set of iso-classes of simple $A$-modules and $A_0 = \oplus_{S \in \mathcal{S}_A}S$. 

For $M\in \mod A$ we denote by $M^k = \oplus_{i=1}^k M$, by $P(M)$ its projective cover, and by $\Omega(M)$ its syzygy. 

For a subcategory $\mathcal{C} \subset \mod A$, we denote by $\findim (\mathcal{C})$ and $\gl (\mathcal{C})$ its finitistic dimension and its global dimension respectively. We use $\add \mathcal{C}$ for the full subcategory of $\mod A$ whose objects are the direct sums of direct summands of the modules from $\mathcal{C}$. 

If $Q = (Q_0,Q_1,\start,\target)$ is a finite connected quiver, $\Bbbk Q$ denotes its associated {\bf path algebra} (we compose paths in $Q$ from left to right). 
Given $\rho$ a path in $\Bbbk Q$, $\length(\rho)$, $\start(\rho)$ and $\target(\rho )$ denote the length, start and target of $\rho$ respectively. If there is a path $\rho$ such that $v= \start(\rho)$ and $w=\target(\rho)$, we say that $w$ is a {\bf successor} of $v$ and $v$ is a {\bf predecessor} of $w$.

For a quiver $Q$, we denote by $J_Q$ the ideal of $\Bbbk Q$ generated by all of its arrows. If $I$ is an admissible ideal of $\Bbbk Q$ ($J^m_Q \subset I \subset J^2_Q$ for some $m \geq 2$), we say that $(Q,I)$ is a {\bf bounded quiver} and the quotient algebra $\frac{\Bbbk Q}{I}$ is the {\bf bound quiver algebra}. A {\bf relation} $\rho$ is a finite sum in $\Bbbk Q$, $\rho = \sum \lambda_i w_i$ where the $\lambda_i$ are scalars (not all zero) and the $w_i$ are paths with $\length(w_i) \geq 2$ such that $\start(w_i) = \start(w_j)$ and $\target(w_i) = \target (w_j)$ if $i \not = j$. 
We recall that an admissible ideal $I$ is always generated by a finite set of relations (for a proof, see Chapter II.2 Corollary 2.9 of \cite{ASS}).

%

For a quiver $Q$, we say that $M = (M^v, T_{\alpha})_{v\in Q_0, \alpha \in Q_1}$ is a {\bf representation} if
\begin{itemize}
\item $M^v$ is a $\Bbbk$-vectorial space for every $v\in Q_0$,
\item $T_{\alpha} :M^{\start(\alpha)} \rightarrow M^{\target(\alpha)}$ is a $\Bbbk$-linear map for every $\alpha \in Q_1$.
\end{itemize}

A representation $(M^v, T_{\alpha})_{v\in Q_0, \alpha \in Q_1}$ is finite dimensional if $M^v$ is finite dimensional for every $v \in Q_0$.

For a path $w = \alpha_1\ldots \alpha_n$ we define $T_{w} = T_{\alpha_n} \ldots T_{\alpha_1}$, and for a relation $\rho = \sum \lambda_i w_i $ we define $T_{\rho} = \sum \lambda_i T_{w_i}$. 

A representation $M = (M^v, T_{\alpha})_{v\in Q_0, \alpha \in Q_1}$ of $Q$ is {\bf bound by} $I$ if we have  $T_{\rho} = 0$ for all relations $\rho \in I$.
 
Let $M= (M^v, T_{\alpha})_{v\in Q_0, \alpha \in Q_1}$ and $\bar{M} = (\bar{M}^v, \bar{T}_{\alpha})_{v\in Q_0, \alpha \in Q_1}$ be two representations of the bounded quiver $(Q, I)$, a {\bf morphism} $f : M \rightarrow \bar{M}$ is a family $f=(f_v)_{v\in Q_0}$ of $\Bbbk$-linear maps $f_v: M^v \rightarrow \bar{M}^v$ such that for all arrow $\alpha$  with $\start(\alpha) = v$, and $\target(\alpha)=w$ the following diagram is commutative:

$$\xymatrix{ M^v \ar[r]^{T_{\alpha}} \ar[d]_{f_v} & M^w \ar[d]^{f_w}  \\ \bar{M}^v \ar[r]^{\bar{T}_{\alpha}}& \bar{M}^w}$$ 

We denote by $\Rep_{\Bbbk}(Q)$ ($\rep_{\Bbbk}(Q)$) the category of (finite dimensional) representations of $Q$. We denote by $\Rep_{\Bbbk} (Q,I)$ ($\rep_{\Bbbk} (Q,I)$) the full subcategory of
$\Rep_{\Bbbk}(Q)$ ($\rep_{\Bbbk}(Q)$) consisting of the representations of
$Q$ bound by $I$.

We recall that there is a $\Bbbk$-linear equivalence of categories $$F: \Rep_{\Bbbk} (Q,I) \rightarrow \Mod \frac{\Bbbk Q}{I}$$ that restricts to an equivalence of categories $G: \rep_{\Bbbk} (Q,I) \rightarrow \mod \frac{\Bbbk Q}{I}$ (for a proof see Chapter III.1, Theorem 1.6 of \cite{ASS}).

For a vertex $v$, we denote by $P_v$ and $S_v$ the indecomposable projective module and the simple module associated to the vertex $v$. For a module $M$, $\supp (M)$ denotes its support, i.e. the vertices $v$ such that $M^v\not=0$.

%
%

\subsection{Igusa-Todorov functions and Igusa-Todorov algebras}

In this section, we exhibit some general facts about the Igusa-Todorov functions for an Artin algebra $A$. We also define the Igusa-Todorov algebras and LIT algebras. The aim is to present the material to be used in the following sections.

\begin{defi}
Let $K_0 (A)$ denote the quotient of the free abelian group generated by one symbol  $[M]$ for each isomorphism class of right finitely generated $A$-module, and relations given by:
\begin{enumerate}
  \item $[M]-[M']-[M'']$ if  $M \cong M' \oplus M''$.
  \item $[P]$ for each projective module $P$.
\end{enumerate}
\end{defi}

Let $\bar{\Omega}: K_0 (A) \rightarrow K_0 (A)$ be the group endomorphism induced by $\Omega$, i.e. $\bar{\Omega}([M]) = [\Omega(M)]$. We denote by $K_i (A) = \bar{\Omega}(K_{i-1}(A))= \ldots = \bar{\Omega}^{i}(K_{0} (A))$ for $i \geq 1$, {\bf the syzygy class groups}. 
We say that a subgroup $G \subset K_0(A)$ is syzygy finite if there is $n \geq 0$ such that $\bar{\Omega}^n(G)$ is finitely generated.  The algebra $A$ is syzygy finite if $K_0(A)$ is syzygy finite.

 For a subcategory $\mathcal{C} \subset \mod A$, we denote by $\langle\mathcal{C}\rangle \subset K_0(A)$ the free abelian group generated by the iso-classes of direct summands of modules of $\mathcal{C}$. In particular, for $M\in \mod A$, $\langle \add M\rangle$ denotes the subgroup of $K_0 (A)$ generated by the iso-classes of indecomposable summands of $M$.

\begin{defi}\label{monomorfismo}(\cite{IT}) 
The \textbf{(right) Igusa-Todorov function} $\phi$ of $M\in \mod A$  is defined as 
\[\phi_{A}(M) = \eta_{\bar{\Omega}}(\langle \add M \rangle).\]
\end{defi}

\begin{prop}(\cite{HLM}, \cite{IT}) \label{it1}
If $M,N\in\mod A$, then we have the following.

\begin{enumerate}
  \item $\phi(M) = \pd (M)$ if $\pd (M) < \infty$.
  \item $\phi(M) = 0$ if $M \in \ind A$ and $\pd(M) = \infty$.
  \item $\phi(M) \leq \phi(M \oplus N)$.
  \item $\phi\left(M^{k}\right) = \phi(M)$ for $k \geq 1$.
  \item $\phi(M) \leq \phi(\Omega(M))+1$.
\end{enumerate}

\end{prop}

\begin{defi} Let $A$ be an Artin algebra, and $\mathcal{C}$ a subcategory of $\mod A$. The $\phi$-dimension of $\mathcal{C}$ is defined as
$$\fidim  (\mathcal{C}) = \sup \{ \phi(M): M \in \mathcal{C} \}.$$
In particular, we denote by $\fidim (A)$ the $\phi$-dimension of $\mod A$.
\end{defi}

We recall the following two theorems from \cite{LM}.

\begin{teo}\label{LM2}(Corollary 3.17. of \cite{LM})
For any Artin algebra $A$ such that $id(A_A) = n < \infty$, it follows
that $$\findim(A) \leq \fidim (A) \leq \psidim (A) \leq n.$$
\end{teo}

\begin{teo}\label{LM1}(Theorem 4.7. of \cite{LM})
For any n-Gorenstein algebra $A$, it follows
that $$\findim(A) = \fidim (A) = \psidim (A) = n.$$
\end{teo}

Now, we recall the notion of LIT-algebra from \cite{BLMV}.

\begin{defi}[\cite{BLMV}]
A subclass $\mathcal{D} \subset \mod A$ is a {\bf $0$-Igusa-Todorov subcategory} if: $\mathcal{D} = \add (\mathcal{D})$, $\Omega (\mathcal{D}) \subset \mathcal{D}$ and $\fidim(\mathcal{D}) = 0$.
\end{defi}

\begin{defi}[\cite{BLMV}]\label{algebra LIT}
An {\bf $(n, V, \mathcal{D})$-Lat-Igusa-Todorov algebra} ({\bf $(n, V, \mathcal{D})$-LIT-algebra}, for short) is an Artin algebra $A$ satisfying the following conditions: n is a non-negative integer, $\mathcal{D} \subset \mod A$ is a $0$-Igusa-Todorov subcategory, and
$V \in \mod A$ such that each $M \in \mod A$ admits a short exact sequence:
\[\xymatrix{0 \ar[r]& V_1 \oplus D_1 \ar[r] & V_0 \oplus D_0 \ar[r] & \Omega^n(M)\ar[r] & 0}\]
such that $V_0, V_1 \in \add A$ and $D_0, D_1 \in \mathcal{D}$.

\end{defi}

\begin{obs}[Definition 2.2 from \cite{W}]
If $\mathcal{D} = \{0\}$ in Definition \ref{algebra LIT}, we say that $A$ is a {\bf n-Igusa-Todorov algebra}.
\end{obs}

The theorem below shows the importance of the previous definitions.

\begin{teo}(Theorem 5.4 of \cite{BLMV}) Let $A$ be an $(n, V, \mathcal{D})$-LIT-algebra. Then
\[\fin (A)  < \infty.\]
\end{teo}

\section{Delooping levels}

\begin{defi}(\cite{G1})
Let $A$ be an Artin algebra and $M \in \mod A$. The delooping level of a $A$-module $M$ is
$$\Dell(M) = \inf \{n\geq 0 : \Omega^n(M) \text{ is a direct summand of } \Omega^{n+1}(N) \text{ in } \underline{\mod} A \text{ for some } N\in \mod A\}.$$
If no such $n$ exist, $\Dell(M) = \infty$. 
\end{defi}

Observe that $\Dell(M) = 0$ if and only if $M \in \Omega(\mod A)$. i.e. $M$ is a torsionless module. 

\begin{defi}\label{dell}
For an Artin algebra $A$, we define its delooping levels as 
\begin{itemize}

\item $\Dell(A) = \sup \{\Dell(S): S \in \mathcal{S}_A\}$ (delooping level (\cite{G1})),

\item $\Dellm(A) = \sup \{\Dell(M): M \in \mod A\}$ (global delooping level).

\end{itemize}
\end{defi}

Observe that $\Dell(A) \leq \Dellm(A)$ for every Artin algebra $A$. And in case $\gl(A) < \infty$, then $\Dell(A) \leq \Dellm(A) < \infty$. \\

We recall the next two results due to V. Gélinas.

\begin{teo}\label{gelinas}(Proposition 1.3 of \cite{G1})
If $R$ is an artinian ring, then
$$ \Fin (R^{op}) \leq \Dell (R).$$
\end{teo}

\begin{teo}\label{gelinasrad^2=0} (Theorem 4.3 of \cite{G2})
Let $A$ be an Artin algebra with radical square zero. Then we have
$$\Fin (A^{op}) = \Dell (A).$$

\end{teo}

As a direct consequence of the first theorem, we have $ \findim (R^{op}) \leq \Dell (R)$ for $R$ an artinian ring. Also, by Theorem 3.1 from \cite{M} we have $\fin(A^{op}) = \Dell(A)$ for radical square zero algebras.

The proposition below shows two properties for the delooping level used in the rest of the article.

\begin{prop} \label{suma}
If $M_1, M_2 \in \mod A$ then:
\begin{enumerate}
\item $\Dell (M_1 \oplus M_2) = \max\{\Dell(M_1), \Dell(M_2)\}$.

\item $\add(M_1) = \add(M_2)$ implies $\Dell (M_1) = \Dell(M_2)$.

\end{enumerate}

\begin{proof}

\begin{enumerate}
\item Assume that $\Dell(M_1) = m_1$ and $\Dell(M_2) = m_2$ with $m_1 \leq m_2$. Then $\Omega^{m_1}(M_1)$ is a direct summand of $\Omega^{m_1+1}(N_1)$ in $\underline{\mod} A$ and $\Omega^{m_2}(M_2)$ is a direct summand of $\Omega^{m_2+1}(N_2)$ in $\underline{\mod} A$, for some $N_1$ and $N_2 \in \mod A$ . 
Clearly $\Omega^{m_2}(M_1)$ is a direct summand of $\Omega^{m_2+1}(N_1)$ in $\underline{\mod} A$, thus $\Omega^{m_2}(M_1 \oplus M_2)$ is a direct summand of $\Omega^{m_2+1}(N_1 \oplus N_2)$ in $\underline{\mod} A$. Therefore $\Dell (M_1 \oplus M_2) \leq \max\{\Dell(M_1), \Dell(M_2)\}$.\\

On the other hand, consider $\Dell(M_1\oplus M_2) = m$. Then $\Omega^{m}(M_1 \oplus M_2)$ is a direct summand of $\Omega^{m+1}(N)$ in $\underline{\mod} A$. This implies that  $\Omega^{m}(M_1)$ and $\Omega^{m}(M_2)$ are direct summand of $\Omega^{m+1}(N)$ in $\underline{\mod} A$. Hence $\max\{\Dell (M_1),\Dell(M_2)\} \leq m = \Dell(M_1\oplus M_2)$, and we are done.  

\item It is a consequence of item 1.
\end{enumerate}
\end{proof}

\end{prop}

As a consequence of the above proposition, we get the following remark.

\begin{obs} For an Artin algebra
$\Dell(A) = \Dell(A_0) \leq \Dellm (A)$
\end{obs}

\begin{prop}
If $A$ is an Artin algebra, $\Dellm(A) = m < \infty$ if and only if $m$ is minimal with the condition $\add(\Omega^m(\mod A)) = \displaystyle \bigcap_{j \geq m} \add ( \Omega^j(\mod A))$.
\begin{proof}
If $\Dellm(A) = m <\infty$, then $\forall M \in \mod A$ and $\forall k \geq 0$ there exists $N_k$ : $\Omega^m(M)$ is a direct summand of $\Omega^{m+k}(N_k)$.
\end{proof}

\end{prop}

The following example shows that the delooping level can be infinite.

\begin{ej}(\cite{KR}) \label{ejmeplo inf} Let $\Bbbk$ be a field, and $q \in \Bbbk$ an element with infinite multiplicative order.
Then $A = A(q)$ is the algebra

$$\frac{\Bbbk \langle x, y, z\rangle} {\langle x^2 , y^2 , z^2 , zy, yx + qxy, zx-xz, yz-xz\rangle}.$$

For each $\alpha \in \Bbbk$, $M (\alpha)$ is the three-dimensional $A$-module with basis
$v, v' , v'' $, where $vx = \alpha v'$ , $vy = v'$ , $vz = v''$ , with $v'$ and $v''$ annihilated
by $x$, $y$ and $z$.

Then the one-point extension algebra $A[M(q)]$ has infinite delooping level (see Corollary 4.2 from \cite{KR}). 

\end{ej}

\subsection{Relations between $\Dell$ and other homological dimensions}

\begin{prop}\label{sizigia finita}
If $A$ is syzygy finite, then $$\Dell(A) \leq \Dellm(A) \leq \fidim(A) < \infty.$$

\begin{proof}
Since $A$ is syzygy finite there is $n \in \mathbb{N}$ such that $\bar{\Omega}:K_{n+k}(A) \rightarrow K_{n+1+k}(A)$ is an isomorphism $\forall k \in \mathbb{N}$. Suppose $n$ is the minimum, then, by Proposition 3.9 from \cite{LMM}, we have $\fidim(A) = n$.

On the other hand, because $A$ is syzygy finite, then $K_n(A)$ is finitely generated. Now, since $\bar{\Omega}:K_{n+k}(A) \rightarrow K_{n+1+k}(A)$ is an isomorphism and $\rk(K_{n+k}) < \infty$  $\forall k \in \mathbb{N}$, for every $M \in \mod A$, there is $N \in \mod A$ such that $\Omega^n(M)$ is a direct summand of $\Omega^{n+1}(N)$. Hence, we conclude that $$\Dellm(A) \leq \fidim (A) = n.$$
\end{proof}
\end{prop}

It was already known that for an Artin algebra, $\findim(A) \leq \fidim(A)$. Now, as a consequence of the previous result, we also have:

\begin{coro}
If $A$ is syzygy finite, then $\findim(A^{op})\leq \fidim(A).$ 

\begin{proof}
Apply  Theorem \ref{gelinas} and Proposition \ref{sizigia finita}.
\end{proof}

\begin{obs}
The inequality $\Dell(A) \leq \fidim (A)$ was first obtained in \cite{GI} (Theorem 1.2) for algebras with $\{\Omega^i(S_j): \forall j\in 1, \ldots, n, \forall i \geq 0 \}  \subset \add T$ for $T \in \mod A$ where $\{S_1, \ldots, S_n\} = \mathcal{S}_A$. 
\end{obs}

\end{coro}

The Proposition \ref{sizigia finita} is not valid in general.

\begin{ej} We exhibit two examples of syzygy infinite algebras. The first one still verify the inequality $\Dell(A) \leq \fidim(A)$ from Proposition \ref{sizigia finita}, and the second one  do not verify such inequality.

\begin{enumerate}

\item Consider Example 5.4 from \cite{BM1}. Let $A = \frac{\Bbbk Q}{I}$ be an algebra where $Q$ is 

$$
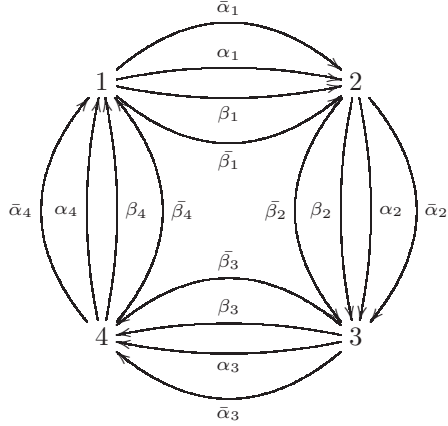
\xymatrix{ 1 \ar@/^8mm/[rrr]^{\bar{\alpha}_1} \ar@/^2mm/[rrr]^{\alpha_1} \ar@/_2mm/[rrr]_{\beta_1} \ar@/_8mm/[rrr]_{\bar{\beta_1}} & &  & 2 \ar@/^8mm/[ddd]^{\bar{\alpha}_2} \ar@/^2mm/[ddd]^{\alpha_2} \ar@/_2mm/[ddd]_{\beta_2} \ar@/_8mm/[ddd]_{\bar{\beta_2}} \\ & &  &\\& & & & \\ 4 \ar@/^8mm/[uuu]^{\bar{\alpha}_4} \ar@/^2mm/[uuu]^{\alpha_4} \ar@/_2mm/[uuu]_{\beta_4} \ar@/_8mm/[uuu]_{\bar{\beta_4}} &  & &  3 \ar@/^8mm/[lll]^{\bar{\alpha}_3} \ar@/^2mm/[lll]^{\alpha_3} \ar@/_2mm/[lll]_{\beta_3} \ar@/_8mm/[lll]_{\bar{\beta_3}}}$$

and $I = \langle \alpha_{i}\alpha_{i+1}-\bar{\alpha}_{i}\bar{\alpha}_{i+1},\ \beta_{i}\beta_{i+1}-\bar{\beta}_{i}\bar{\beta}_{i+1},\ \alpha_{i}\bar{\alpha}_{i+1},\ \bar{\alpha}_{i}\alpha_{i+1},\ \beta_{i}\bar{\beta}_{i+1},\ \bar{\beta}_{i}\beta_{i+1},\text{ for } i \in \mathbb{Z}_4, \ J^3 \rangle$

We know that $\fidim(A) = \fidim(A^{op}) =  \infty$. On the other hand it is clear that $S_i \subset P_{i+2}$ for $i \in \mathbb{Z}_4$, hence $\Dell(A) < \fidim(A) = \fidim(A^{op})$.

\item Consider Example 4.2 from \cite{BM}. Let $C_{p,q} =\frac{\Bbbk Q}{I_{p,q}}$ be the finite dimensional algebra where $Q$ is the following quiver

$$
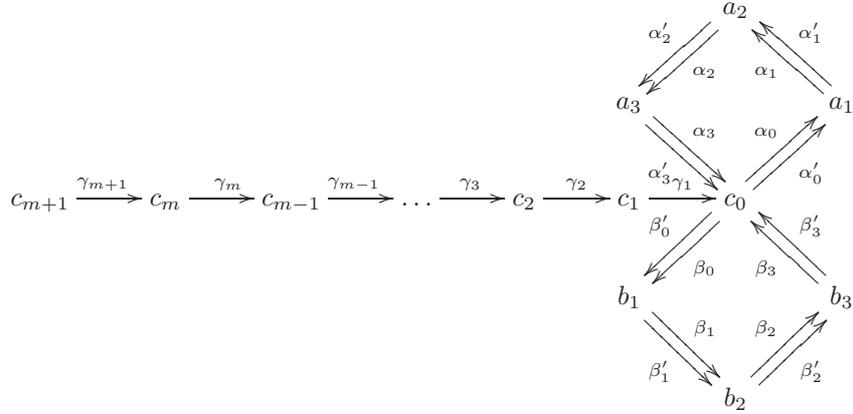
\xymatrix{ &&&&&& a_2 \ar@<.5ex>[dl]^{\alpha_2} \ar@<-.5ex>[dl]_{\alpha'_2} && \\&&&&& a_3 \ar@<.5ex>[dr]^{\alpha_3} \ar@<-.5ex>[dr]_{\alpha'_3} && a_1 \ar@<.5ex>[ul]^{\alpha_1} \ar@<-.5ex>[ul]_{\alpha'_1}& \\ c_{m+1} \ar[r]^{\gamma_{m+1}} &  c_m \ar[r]^{\gamma_m}&c_{m-1}\ar[r]^{\gamma_{m-1}}& \ldots \ar[r]^{\gamma_3}& c_2 \ar[r]^{\gamma_2}& c_1 \ar[r]^{\gamma_1}& c_0\ar@<.5ex>[ur]^{\alpha_0} \ar@<-.5ex>[ur]_{\alpha'_0} \ar@<.5ex>[dl]^{\beta_0} \ar@<-.5ex>[dl]_{\beta'_0} &&  \\&&&&& b_1 \ar@<.5ex>[dr]^{\beta_1} \ar@<-.5ex>[dr]_{\beta'_1} && b_3 \ar@<.5ex>[ul]^{\beta_3} \ar@<-.5ex>[ul]_{\beta'_3}& \\&&&&&& b_2 \ar@<.5ex>[ur]^{\beta_2} \ar@<-.5ex>[ur]_{\beta'_2} && }$$

and $I_{p,q}$ is generated by

\begin{itemize}
\item $\gamma_{i+1}\gamma_{i}\ \forall i=1,\ldots,m$,

\item $\alpha_{i}\alpha_{i+1}, \beta_{i}\beta_{i+1}, \alpha'_{i}\alpha'_{i+1}, \beta'_{i}\beta'_{i+1}\ \forall i \in \mathbb{Z}_4$,

\item $\alpha'_{i}\alpha_{i+1} - \alpha_{i}\alpha'_{i+1}, \beta'_{i}\beta_{i+1}-\beta_{i}\beta'_{i+1}\ \forall i \in \mathbb{Z}_4$,

\item $ \alpha_3\beta_0 , \alpha'_3\beta'_0, \beta_3\alpha_0, \beta'_3\alpha'_0$ and
\item $\alpha'_3\beta_0 -\alpha_3\beta'_0, \beta'_3\alpha_0-p\beta_3\alpha'_0, \gamma_1\alpha_0-q\gamma_1\alpha'_0, \gamma_1\beta'_0-\gamma_1\beta_0$ with $\mathbb{Q}[p,q] \cong \mathbb{Q}[x,y]$. 
\end{itemize}

We know that $\fidim(C_{p,q}) = 5$. 

On the other hand, we have

\begin{itemize}

\item $\Dell(S_{c_i}) = 0$ for $i = 1, \ldots m$

\item $\Dell(S_{c_0}) =  \Dell(S_{a_i}) = \Dell(S_{b_j}) = 0$ for $i=1,2,3$, $j=1,2,3$ (use the calculations from Example 4.2 of \cite{BM}).
\item $\Dell(S_{c_{m+1}}) = m+2$ (use the calculations from Example 4.2 of \cite{BM}).

\end{itemize}
 
Hence $5 = \fidim(C_{p,q}) < \Dell(C_{p,q}) = m+2$. Hence the difference $\Dell(C_{p,q}) - \fidim (C_{p,q})$ can be arbitrarily large.

\end{enumerate}

\end{ej}

Igusa-Todorov algebras and LIT algebras were introduced as tools to prove the finitistic dimension conjecture, for large families of algebras.
The following example shows that both techniques encompass different families of algebras, none of which is included in each other.     
 
\begin{ej} Note that

\begin{enumerate}
\item The opposite algebra of $A(q)$ from Example \ref{ejmeplo inf} is a radical cube zero algebra. Hence it is an Igusa-Todorov algebra (Theorem 1.3. from \cite{W}).

\item Consider the algebra $A = \frac{\Bbbk Q}{I}$ of Example 4 from \cite{BM2}, where $Q$ is the following quiver
$$
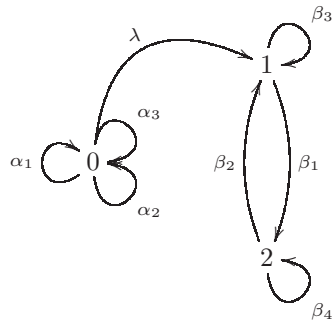
\xymatrix{ & & 1 \ar@(u,r)^{\beta_3} \ar@/^3mm/[dd]^{\beta_1} \\ 0 \ar@(dl,lu)^{\alpha_1} \ar@(d,r)_{\alpha_2} \ar@(u,r)^{\alpha_3} \ar@/^10mm/[urr]^{\lambda} & & \\ & & 2 \ar@(d,r)_{\beta_4}  \ar@/^3mm/[uu]^{\beta_2}  }$$ 
and $I = \langle \alpha_i \alpha_j + \alpha_j \alpha_i\ \forall i=1,2,3\text{ and }\forall j=1,2,3, \alpha_i \lambda \ \forall i=1,2,3, \lambda \beta_1, \lambda \beta_3, \beta_i\beta_j \ \forall i=1,2,3,4 \text{ and } \forall j=1,2,3,4 \rangle$. 

We know that $A$ is not a LIT algebra (Example 4 from \cite{BM2}). On the other hand, $\Dell(A) = \Dell(A^{op}) = 0$ since every simple module is a submodule of a projective module.

\end{enumerate}

For references, see \cite{BLMV}.

\end{ej}

\subsection{Dell on Gorenstein algebras}

The next lemma gives a new way to compute the minimum $n$ for which a Gorenstein algebra is $n$-Gorenstein, using the Igusa-Todorov functions.

\begin{lema} \label{Lemita}
If $A$ is a Gorenstein algebra, then $A$ is $n$-Gorenstein (with $n$ minimal) if and only if $K_n(A)$ is the first level of the syzygy class groups generated by Gorenstein projective modules.

\begin{proof}
($\Rightarrow$) We know that if $A$ is $n$-Gorenstein (with $n$ minimal), then $K_n(A)$ is generated by classes of Gorenstein projective modules and $\fidim(A) = n$ (see Theorem 4.7 from \cite{LM}). If $K_{n-1}(A)$ is generated by Gorenstein projective modules, then by item 5 from Proposition \ref{it1} and Corollary 4.1 from \cite{LM} $\fidim (A) \leq n-1$, which is absurd.

($\Leftarrow$) If $K_n(A)$ is generated by classes of Gorenstein projective modules and $K_{n-1}(A)$ it is not, then $\fidim (A) = n$, thus $A$ is $n$-Gorenstein (with $n$ minimal).
\end{proof}

\end{lema}

\begin{teo}\label{n-Goreinstein}
If $A$ is $n$-Gorenstein then $\Dellm(A) = n$.
\begin{proof}

Since $A$ is $n$-Gorenstein, $\Omega^{n}(M)$ is a Gorenstein projective module for all $M \in \mod A$. We know that for every Gorenstein projective module $G$ there exists another Gorenstein projective module $\bar{G}$ such that $\Omega(\bar{G}) = G$, thus $\Dell(M) \leq n$.\\

Suppose $\Dellm(A) = m<n$. Then for every $M \in \mod A$, there exists $N$ such that $\Omega^m(M)$ is direct summand of the module $\Omega^{m+1}(N)$ (in $\underline{\mod} A$). Since $\Omega^{n-m-1}( \Omega^m(M))=\Omega^{n-1}(M)$ is direct summand of $\Omega^{n-m-1} (\Omega^{m+1}(N)) = \Omega^n(N)$, then $\Omega^{n-1}(M)$ is Gorenstein projective for every $M \in \mod A$, a contradiction with Lemma \ref{Lemita}.
\end{proof}

\end{teo}

We recall that if $M \in \mod A$ is a direct summand of $\Omega(\bar{N})$ for some $\bar{N} \in \mod A$, there exists $N \in \mod A$ such that $\Omega(N) = M$. \\

\begin{coro}\label{Selfinjective}
 If $A$ is an Artin algebra, the following statements are equivalent

\begin{itemize}
\item[1.]  $A$ is selfinjective.

\item[2.] $\Dell(DA_A) = 0$.

\item[3.] $\Dellm(A) = 0$.
\end{itemize} 
  
\begin{proof}

First note that if $\Dell(I) = 0$ for an injective $A$-module $I$, then there is a short exact sequence 
$$\delta: 0\rightarrow I \rightarrow P \rightarrow M \rightarrow 0$$with $P$ a projective module. Hence $I$ is projective since $\delta$ splits.\\

\noindent If $\Dellm(A) = 0$, then  $\Dell(DA_A) = 0$, therefore, by previous remark, $DA_A$ is a projective module.  This implies that $A$ is selfinjective; in other words, $A$ is $0$-Gorenstein. Therefore, by Theorem \ref{n-Goreinstein}, we get that $\Dellm(A) = 0$, proving the three equivalences.


\end{proof}

\end{coro}

\begin{prop}\label{id finita}
Let $A$ be an Artin algebra with $\id(A)< \infty $, then $$\findim(A)\leq \fidim (A) \leq \psidim (A)\leq  \id(A) \leq \Dell(A) \leq \Dellm(A).$$

\begin{proof}
Since $\pd_{A^{op}}(D(A_A)) = 
\id_A (A_A)<\infty $, by Theorem \ref{gelinas}, we have
$$\pd_{A^{op}}(D(A_A)) \leq \findim(A^{op})\leq \Dell(A).$$
On the other hand, by Theorem \ref{LM2}, we have 
$$\findim(A) \leq \fidim(A)\leq \psidim(A) \leq \pd_{A^{op}}(D(A_A)).$$
Therefore, gluing both chains of inequalities:
$$\findim(A) \leq \fidim(A)\leq \psidim(A) \leq \pd_{A^{op}}(D(A_A))  \leq \findim(A^{op})\leq \Dell(A).$$\end{proof}
\end{prop}

\noindent The following result was mentioned as Example 1.18 of \cite{G1}.
\begin{coro}
If $A$ is $n$-Gorenstein then $\Dell(A) = n$.
\begin{proof}
By Theorem \ref{n-Goreinstein}, $\Dell(A) \leq \Dellm(A) = n$.
On the other hand $n = \pd(D(A_A)) \leq \findim(A^{op})\leq \Dell(A)$, by Theorem \ref{gelinas}.
\end{proof} 

\end{coro}

Since $A$ is a $n$-Gorenstein algebra if and only if $A^{op}$ is a $n$-Gorenstein algebra, we deduce that for Gorenstein algebras $\Dell(A) =\Dellm(A) = \Dellm(A^{op}) =  \Dell(A^{op})$. This remark motivates the following example.

\begin{ej} \label{simetricas}
Are $\Dell(A)$ or $\Dellm(A)$ symmetric homological measures?

\begin{itemize}
\item It is known that if $A = \frac{\Bbbk Q}{J^2}$, then $\fin (A) \not = \fin (A^{op})$, in general. Since $\Dell(A) = \fin (A^{op})$ by Theorem \ref{gelinasrad^2=0}, then the delooping level is not a symmetric homological measure. For instance, consider $Q$ the following quiver
$$ \xymatrix{ 1 \ar[r]  \ar@(dl,lu) & 2 \ar[r] & 3 \ar[r] & \ldots \ar[r] & n-1 \ar[r] & n \ar[r] & n+1 & }.$$

It is clear that $\findim(A) = n$, hence $\Dell(A^{op})= \Dellm(A^{op})=n$.
On the other hand, since there are no $A^{op}$-modules with positive projective dimension, then $\Dell(A) = 0$ and by Corollary \ref{radical^2=0 Dell} $\Dellm(A) = 1$.  

\item $A[M(q)]^{op} $ from Example \ref{ejmeplo inf} verifies that $\Dell(A[M(q)]^{op}) = 0$ since every simple module is a submodule of $A[M(q)]^{op}$. It shows that the difference between right and left delooping levels can be infinite. 
\end{itemize}

\end{ej}

\section{Monomial algebras}

In this section, we study the delooping levels in monomial algebras, with a special emphasis in truncated path algebras. We prove a generalization of Theorem \ref{gelinasrad^2=0} for truncated path algebras. Also, we show an example of a family $\{A_n\}_{n\in \mathbb{N}}$ of monomial algebras where the difference $\Dell (A_n)-\Fin(A_n)$ can be arbitrary large. 

\subsection{Radical square zero algebras}

We recall that V. Gélinas prove in \cite{G2} that $ \Dell(A) = \ \Fin(A^{op}) = \findim(A^{op})$ for a radical square zero algebra (Theorem \ref{gelinasrad^2=0}). Now, we show a similar result for the global delooping level $(\Dellm)$.

\begin{obs}\label{sizigia radical^2=0}
Let $A$ be  a radical square zero algebra. 
\begin{itemize}
\item If $M \in \mod A$, then $\Omega(M)$ is a direct summand of $\Omega(\topp (M))$. 
\item If $M$ is an indecomposable non-simple non-projective module, then $\Dell(M) > 0$.

\item If $v$ is a source, then $\Dell(S_v) > 0$.
\end{itemize}

\end{obs}

\begin{prop}\label{radical^2=0}
Let $A$ be a radical square zero algebra. If $M$ is an indecomposable non-simple non-projective module, then 
$$ 1 \leq \Dell(M) \leq \max\{1, \Dell(\topp(M))\}.$$

\begin{proof}
First observe that, by Remark \ref{sizigia radical^2=0}, $\Omega^n(M)$ is a direct summand of $\Omega^n(\topp(M))$ and $\Dell(M) > 0$.

If $\Dell(\topp(M)) = 0$, then $\topp(M) = \Omega(N)$ for some $N \in \mod A$.  Since $\Omega(M)$ is a direct summand of $\Omega(\topp(M))$, we get that $\Omega(M)$ is a direct summand of $\Omega^2(N)$. Hence $\Dell(M) \leq 1$. On the other hand, by  Remark \ref{sizigia radical^2=0} item 2, $\Dell(M) > 0$.

If $\Dell(\topp(M)) = n \geq 1$, then there is a $A$-module $N$ such that $\Omega^n(\topp(M))$ is a direct summand of $\Omega^{n+1}(N)$. By Remark \ref{sizigia radical^2=0} item 1, $\Omega^n(M)$ is a direct summand of $\Omega^{n+1}(N)$. Hence, $\Dell(M) \leq n$.

\end{proof}

\end{prop}

Observe that by Proposition \ref{suma}  and Proposition \ref{radical^2=0}, $\Dell(M) > 0$ if there exists an indecomposable direct summand of $M$ that is not a projective or simple. 

By Corollary \ref{Selfinjective}, Remark \ref{sizigia radical^2=0} and Proposition \ref{radical^2=0}, we get the next result follows.

\begin{coro}\label{radical^2=0 Dell}
Let $A$ be a non-selfinjective radical square zero algebra.

\begin{itemize}
\item If $\Dell(A) = 0 $, then $\Dellm(A) = 1$.
\item If $\Dell(A) > 0$, then $\Dell(A) = \Dellm(A)$.

\end{itemize}

\end{coro}

\subsection{Truncated path algebras}

This subsection is devoted to extend Gelinas' Theorem 4.3 from \cite{G2} (Theorem \ref{gelinasrad^2=0}) to truncated path algebras (Theorem \ref{TeoTruncada} and Corollary \ref{Dell=fin}).

%
%
%
%
%
%
%

\begin{defi}

Let $A = \frac{\Bbbk Q}{I}$. Given a non-trivial path $\rho \not \in I$, for a non-trivial path $\mu_{\rho} \not \in I$, we say that it is a {\bf right complementary path} of $\rho$ if it verifies

\begin{itemize}

\item $\rho \mu_{\rho} \in I$ or $\tg(\mu_{\rho})$ is a sink, and

\item $\rho \mu \not \in I$ $\forall \mu$ such that there is a non-trivial path $\mu'$ with $\mu \mu' = \mu_{\rho}$.
\end{itemize}

\end{defi}

Analogously, we define a {\bf left complementary path} of $\rho$.

\begin{obs}\label{sizigia uniseriales}
Let $A = \frac{\Bbbk Q}{J^l}$ be a truncated path algebra. We recall from Lemma 5.11 of \cite{BH-ZT} that $\Omega(M) = \oplus_{\mu \in I} \mu A$ for every $A$-module $M$. In particular, for ${\rho}A$ we have $$ \Omega({\rho}A) = \bigoplus_{\mu \in I_{\rho}} \mu A,$$ where $I_{\rho} = \{ \mu: \mu \text{ is a right complementary path of } \rho \}$.

\end{obs}

\begin{defi}
Given a path algebra $A = \frac{\Bbbk Q}{I}$, we say that $\tau = \rho_0, \rho_1, \ldots, \rho_k$, a sequence of paths in $\Bbbk Q$, is a {\bf trajectory} from $\rho_0$ to $\rho_k$ if

\begin{itemize}
\item $\rho_{i+1}$ is a right complementary path of $\rho_i$ $\forall i = 0, 1, \ldots, k-1$, 

\item $\rho_{i}$ is a left complementary path of $\rho_{i+1}$ $\forall i = 0, 1, \ldots, k-1$. 
\end{itemize}

We also say that the trajectory $\tau =  \rho_0, \rho_1, \ldots, \rho_k$ has {\bf length} $k$

For $k \in \mathbb{Z^{+}}$ we denote by 
\begin{itemize}
\item $Tr^k(\rho) = \{\tau: \tau  \text{ is a trajectory of length } k \text{ starting at } \rho  \}$,
\item $Tr_k(\rho) = \{\tau: \tau \text{ is a trajectory of length } k \text{ ending at } \rho   \}$.
\end{itemize}
\end{defi}

The next corollary follows directly from Remark \ref{sizigia uniseriales}.

\begin{coro} \label{trajectory}
Let $A = \frac{\Bbbk Q}{J^l}$ be a truncated path algebra. Then $$ \Omega^k(\rho A) = \bigoplus_{\tau \in Tr^k(\rho)} \mu_{\tau} A$$
where $\mu_{\tau}$ is the last path of $\tau$.
\end{coro}

Since being truncated is a symmetric concept, i.e. if $A$ is a truncated path algebra, then $A^{op}$ is also a truncated path algebra, we have the following remark.

\begin{obs} \label{trajectoria dual}
Let $A = \frac{\Bbbk Q}{J^l}$ be a truncated path algebra. Given a trajectory $\tau = \rho_0, \rho_1, \ldots, \rho_k $, then $\tilde{\tau} = \tilde{\rho_k}, \tilde{\rho_{k-1}}, \ldots, \tilde{\rho_0}$ is a trajectory in $A^{op}$, where $\tilde{\rho_i}$ is the opposite path of $\rho_i$ but in $A^{op}$ $\forall i = 0,1, \ldots k$.
\end{obs}

\begin{prop}\label{ni proy ni iny}
Let $A = \frac{\Bbbk Q}{I}$ be a path algebra. Given a trajectory $\tau = \rho_0, \rho_1, \ldots, \rho_k$ with $k \geq 2$, then $\rho_i\rho_{i+1} = 0$ $\forall i = 0,\ldots, k-1$. As a consequence

\begin{itemize}
\item $\rho_i A$ is not projective $\forall i = 0, \ldots, k-1$, and

\item  $\rho_i A$ is not injective $\forall i = 1, \ldots, k$.
\end{itemize} 

\begin{proof}
The first part of the result follows from the fact that $\rho_i \rho_{i+1} = 0$ $\forall i= 0, \ldots, k-2$, since $\rho_{i+1}$ is a right complementary of $\rho_i$ and $\rho_{k-1} \rho_k = 0$ since $\rho_{k-1}$ is a left complementary of $\rho_k$. 

The second one follows analogously using Remark \ref{trajectoria dual}.
\end{proof}

\end{prop}

We can also write the previous corollary in the following way. If there is a trajectory $\tau = \rho_0, \rho_1, \ldots, \rho_k$, then $\rho_k A$ is a direct summand of $\Omega^{k}(\rho_0A)$.

\begin{coro} \label{pd trajectoria}
Let $A = \frac{\Bbbk Q}{J^l}$ be a truncated path algebra. Given a trajectory $\tau = \rho_0, \rho_1, \ldots, \rho_k$ with $k \geq 2$, then $\pd(\rho_0A)\geq k$.

\begin{proof}
The result follows from Proposition \ref{ni proy ni iny} and Corollary \ref{trajectory}.
\end{proof}

\end{coro}

\begin{defi}
Let $A = \frac{\Bbbk Q}{J^l}$ be a truncated path algebra. Given an $A$-module $M$, we define $$\chi(M) = \min\{k : \Omega^k(M) = (\bigoplus_{i\in I} \rho_i A )\oplus P \text{ where } P\in \mathcal{P}_A, \rho_iA \not \in \mathcal{P}_A \text{, and } \ \target(\rho_i) \text{ is a succesor of a cycle}\ \forall i \in I \}.$$

\end{defi}

The following proposition summarizes some elementary properties of the map $\chi$.

\begin{prop} \label{chi sizigia} \label{cota con el chi}
Let $A = \frac{\Bbbk Q}{J^l}$ be a truncated path algebra, then

\begin{enumerate}
\item $\chi(M)$ is always a natural number.  

\item If $M = \oplus_{i\in I} M_i$, then $\chi(M) = \max_{i \in I}\{\chi(M_i)\}$.

\item If $\tg(\rho)$ is a successor of a cycle $Q$, then there a is path $\bar{\rho}$ such that $\rho A = \bar{\rho} A$ and $Tr_n(\bar{\rho}) \not =  \emptyset$, $\forall n \in \mathbb{N}$. 

\item $\Omega^{\chi(M)}(M) \in \add (\Omega^k(\mod A))$ $\forall k \in \mathbb{N}$.

\item Given $M \in \mod A$, then $$\Dell(M) \leq \chi(M) \leq \chi(\Omega(M))+1.$$

\end{enumerate}
\begin{proof}

$1.$, $2.$, and $3.$ follow directly from the definitions. $4.$ is a direct consequence of $3.$. Finally, $5.$ follows from $4.$, Proposition \ref{suma}, and Corollary \ref{trajectory}.
\end{proof}

\end{prop}

\begin{prop} \label{chi}
Let $A = \frac{\Bbbk Q}{J^l}$ be a truncated path algebra and $M \in \mod A$. Then there is a path $\mu_i$ such that
$$\chi(M)-1 \leq \pd(\tilde{\mu_i} A^{op}) < \infty. $$

\begin{proof}

We have three different cases.

\begin{itemize}
\item First consider $\chi(M) \leq 1$. This case is clear.

\item If $\chi(M) \geq 3$, there is a path $\mu_{\chi(M)-2}$ such that $\mu_{\chi(M)-2}A$ is non-projective summand of $\Omega^{\chi(M)-1}{M}$ and $\target(\mu_{\chi(M)-2})$ is not a succesor of any cycle. It is easy to see that there is a trajectory $\tau = \mu_0, \mu_1, \ldots, \mu_{\chi(M)-2}, \mu_{\chi(M)-1}$. Since $\target(\mu_{\chi(M)-2})$ is not a succesor of any cycle, $\tilde{\mu}_{\chi(M)-1}A^{op}$ has finite projective dimension, and by Corollary \ref{pd trajectoria} $\pd(\tilde{\mu}_{\chi(M)-1} A^{op}) \geq \chi(M)-1$.

\item Finally, if $\chi(M) = 2$, we present a particular proof since we cannot use Corollary \ref{pd trajectoria}. In this case there is a path $\mu_0$ such that $\mu_0A$ is a non-projective summand of $\Omega(M)$ and $\target(\mu_0)$ is not a successor of any cycle. Since $\mu_0A$ is not projective then there exists a trajectory $\mu_0, \mu_1$ with $\mu_0 \mu_1 = 0$. Observe that $\target(\mu_0) = \start(\mu_1)$, therefore $\tilde{\mu}_1 A^{op}$ has finite projective dimension and it is not projective since $\mu_0 \mu_1 = 0$.

\end{itemize}

%
%

\end{proof}
\end{prop}

\begin{nota}
Given a non-zero path $\rho = \alpha_1 \ldots \alpha_n$ in $Q$ with $\alpha_i \in Q_1$ $\forall i =1 \ldots n$, we denote by $M_{\rho} = M_{\alpha_1,\ldots, \alpha_n} = (V^v, T_{\alpha}
)$, where 

\begin{itemize}
\item $\{e_{i}: \start(\alpha_i) = v\}$ is a basis of $V^v$ for $v\not = \target(\alpha_n)$, 

\item $\{e_{i}: \start(\alpha_i) = v\} \cup \{e_{n+1}\}$ is a basis of $V^n$ for $v = \target(\alpha_n)$.

\item $T_{\alpha} (e_{i}) = \left\{\begin{array}{ll} e_{i+1} & \text{ if } \alpha = \alpha_i, \\ 0 & \text{ otherwise. } \end{array} \right. $
\end{itemize}

\end{nota}

Observe that the modules $M_{\rho}$ are uniserial.

\begin{obs} \label{uniserial}
Let $A = \frac{\Bbbk Q}{J^l}$ be a truncated path algebra. Consider $\rho = \alpha_1 \ldots \alpha_n$ a non-trivial path in $Q$,  where  $\alpha_1, \ldots, \ \alpha_n$ are arrows. Then $\rho A$ is a direct summand of 
$$\left\{ \begin{array}{ll}
\Omega(M_{\alpha_1 \ldots \alpha_{n-1}}) & \text{ if } n \geq 2, \\ \Omega(S_{\start(\rho)}) & \text{ if } n = 1. \end{array} \right.$$

\end{obs}

\begin{coro} \label{mayor chi}
Let $A = \frac{\Bbbk Q}{J^l}$ be a truncated path algebra with $\findim (A^{op}) \geq 1$. Given $M \in \mod A$, then there exists a module $N \in  \mod A^{op}$ such that $$\chi(M) \leq \pd_{A^{op}}(N) < \infty.$$

\begin{proof}
If $\chi(M) = 1$, there is nothing to prove. If $\chi(M) \geq 2$, by the proof of Proposition \ref{chi}, we know that there is a path $\mu_i$ such that $\chi(M) \leq \pd(\tilde{\mu}_i A^{op})+1 < \infty$. Now using Remark \ref{uniserial}, there exist $N', N'' \in \mod A^{op}$ with $\Omega(N') = \tilde{\mu}_iA^{op} \oplus N''$. Because $\Omega(\mod A^{op})$ is closed by direct summands, therefore there exists $N \in \mod A^{op}$ such that $\Omega(N) = \tilde{\mu}_iA^{op} $.  
\end{proof}

\end{coro}

\begin{prop} \label{cota uniserial}
Let $A = \frac{\Bbbk Q}{J^l}$ be a truncated path algebra. If $M \in \mod A$, then $$\Dell(M) \leq \max \left\lbrace \Dell\left(\bigoplus_{\rho \text{: } \length(\rho)\geq 1} M_{\rho}\oplus A_0\right), 1\right\rbrace.$$

\begin{proof}

Since $\displaystyle \Omega(M) = \bigoplus_{\rho \in R} \rho A$, $\forall M \in \mod A$ (Lemma 5.11 from \cite{BH-ZT}), then $ \Omega(M)$ is a direct summand of $\displaystyle \Omega(\bigoplus_{\mu \in R'} M_{\mu}\oplus A_0)$. The result follows.
\end{proof}
\end{prop}

Now, we can provide proof of this article's main result.

\begin{teo}\label{TeoTruncada}
Let $A = \frac{\Bbbk Q}{J^l}$ be a truncated path algebra.
If $\findim(A^{op}) > 0$, then $$\Dell(A) = \Dellm(A) = \findim(A^{op}).$$

\begin{proof}
Since $\findim(A^{op}) > 0$, then by Corollary \ref{mayor chi} and Proposition \ref{chi sizigia} for all $ M \in mod A$, there exists $N \in \mod A^{op}$ such that
$$\Dell(M) \leq \chi (M) \leq \pd_{A^{op}}(N) < \infty$$
This implies $\Dell(M) \leq \findim(A^{op})$ for all $M \in \mod A$. Hence $\Dellm(A) \leq \findim(A^{op}) \leq \Dell(A)$, where the second inequality is given by Theorem \ref{gelinas}.

On the other hand, we recall that $\Dell(A) \leq \Dellm(A)$.

\end{proof}
\end{teo}

The following corollary is straightforward from the previous theorem.

\begin{coro}\label{S->M}
Let $A = \frac{\Bbbk Q}{J^l}$ be a truncated path algebra with $\findim(A^{op}) > 0$. For all $M \in \mod A$, we have $$\Dell(M) \leq \Dell(\oplus_{v\in Q_0} S_v).$$
\end{coro}

As a consequence of the previous theorem and Theorem 2.2 from \cite{BH-ZT}, we have:

\begin{coro} \label{Dell=fin}
Let $A = \frac{\Bbbk Q}{J^l}$ be a truncated path algebra, then $\Dell(A) = \Fin(A^{op})$.
\end{coro}

Recall that if $A = \frac{\Bbbk Q}{J^l}$ is a truncated path algebra, then $\findim(A^{op}) = 0$ if and only if $Q$ has no sources. Hence every vertex $v$ is a successor of a cycle, $\chi(S_v) = 0$ and $\Dell(A) = 0$.

\begin{prop}

If $A = \frac{\Bbbk Q}{J^l}$ is a truncated path algebra with $\findim(A^{op}) = 0$, then 

\begin{enumerate}
\item If $A$ is selfinjective, then $\Dellm(A) = 0$.

\item If $A$ is not selfinjective, then $\Dellm(A) = 1$.

\end{enumerate}

\begin{proof}

Let $A = \frac{\Bbbk Q}{J^l}$.

\begin{enumerate}
\item If $A$ is selfinjective, we know that $\Dellm(A) = 0$ (Corollary \ref{Selfinjective}).

\item If $A$ is not selfinjective, we know that $\Dellm(A) > 0$ (Corollary \ref{Selfinjective}). Now, since $\findim(A^{op}) = 0$ there are no sources in $Q$. Hence $\chi (M) \leq 1$ for every $M \in \mod A$. By Proposition \ref{cota con el chi} $\Dellm(A) = 1$.
\end{enumerate}
\end{proof}
\end{prop}

Note that the previous result is not valid in general, recall Example \ref{simetricas}. The second item from that example gives an algebra $A(q)$ with $\Dell(A(q)) = 0$ but $\Dell(M(q)) = \infty$, hence $\Dellm(A) = \infty$.

\subsection{Example}
We recall that in \cite{H-Z} B. Huisgen-Zimmermann proved that $\findim (A) \leq \Fin(A) \leq \findim(A) + 1$ for monomial algebras (Corollary B), and she also showed an example of a monomial algebra $A$ where $\Fin(A) = \findim (A)+1$. Barei, Goodearl, and Huisgen-Zimmermann showed an example of an algebra $A$ where $\Dell(A)-\Fin(A) = 1$. This suffices to show that such differences may be arbitrarily large. This is due the following observation of Rickard: On tensor products, the finitistic dimensions and the delooping levels behave additively. This family of algebras are not monomial in general.

In this section, we exhibit an example of a family of algebras which shows that Theorem \ref{TeoTruncada} cannot be generalized to monomial algebras. The example also shows that $\findim(A_n^{op}) \leq \lfloor \frac{k}{4} \rfloor +1 $ and $\Dell(A_n) \geq \lfloor \frac{n}{4} \rfloor$. Hence we can choose $k$ and $n$ so the difference $\Dell (A_n)-\findim(A_n^{op})$ (and also $\Dell (A_n)-\Fin(A_n^{op})$) can be arbitrary large also for monomial algebras. 

\begin{ej}\label{contraej}
Consider the algebra $A_n=\frac{\Bbbk Q_n}{I_n}$, where $Q_n$ is the following quiver
$$
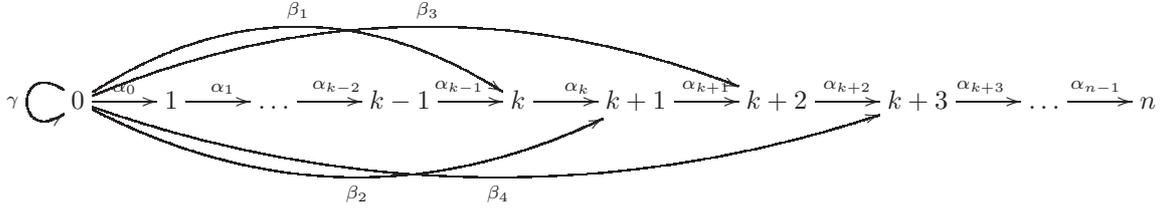
\xymatrix{ &  & &  &  &  &  & &  & &   \\ & 0 \ar@(ul,dl)_{\gamma} \ar[r]^{\alpha_0} \ar@/^10mm/[rrrr]^{\beta_1} \ar@/_10mm/[rrrrr]_{\beta_2} \ar@/^10mm/[rrrrrr]^{\beta_3} \ar@/_10mm/[rrrrrrr]_{\beta_4} & 1 \ar[r]^{\alpha_1}  &\ldots \ar[r]^{\alpha_{k-2}} & k-1 \ar[r]^{\alpha_{k-1}} & k\ar[r]^{\alpha_k} & k+1 \ar[r]^{\alpha_{k+1}} & k+2 \ar[r]^{\alpha_{k+2}} & k+3 \ar[r]^{\alpha_{k+3}} & \ldots \ar[r]^{\alpha_{n-1}}  & n   \\ }$$
and the ideal $I_n = \{\alpha_i\alpha_{i+1}\alpha_{i+2}\alpha_{i+3} \ \forall i =1, \ldots n-4, \beta_i \alpha_{k-1+i}, \gamma \beta_i \ \forall i = 1,2,3,4, \gamma^2,  \alpha_0 \alpha_1 \alpha_2 \}$.\\ 

Since $S_2 \subset P_0$, and $S_i \subset P_{i-3}\ \forall i =4,5, \ldots n-3$, then $\Dell(S_2) = \Dell(S_4) = \Dell(S_5) = \ldots = \Dell(S_{n-1}) = \Dell(S_{n}) = 0$. Thus we have to compute $\Dell(S_0), \Dell(S_1), \Dell(S_3)$ to know $\Dell(A_n)$.

First note that the list of indecomposable modules of $\Omega(\mod A_n)$ is the following
$M_{\alpha_0\alpha_1}, M_{\alpha_1}$, $S_i$, $M_{\alpha_{i-1}}$, $M_{\alpha_{i-2} \alpha_{i-1}}$ $\forall i \geq 4$.

Note also that: 

\begin{itemize}

\item $\Omega(S_0) = M_{\alpha_1} \oplus M_{\alpha_0\alpha_1} \oplus  (\oplus_{i=0}^3 S_{k+i})$

\item $\Omega(M_{\alpha_0\alpha_1}) = M_{\alpha_0\alpha_1} \oplus (\oplus_{i=0}^3 S_{k+i})$,

\item $\Omega(M_{\alpha_1}) = M_{\alpha_3}$,

\item $\Omega(S_i) = \left\{ \begin{array}{cc} M_{\alpha_{i+1}\alpha_{i+2}} & \text{ if } i+3 < n, \\ P & \text{ if } i+3 \geq n. \end{array}\right.$

\item $\Omega(M_{\alpha_{i-1}}) = \left\{ \begin{array}{cc} M_{\alpha_{i+1}} & \text{ if } i+2 < n, \\ P & \text{ if } i+2 \geq n. \end{array}\right.$

\item $\Omega(M_{\alpha_{i-2} \alpha_{i-1}}) = \left\{ \begin{array}{cc} S_{i+1} & \text{ if } i+1 < n, \\ P & \text{ if } i+1 \geq n. \end{array}\right.$

\end{itemize}

Where $P$ is a projective module.

Consider $\mathcal{T}_0 = \add(M_{\alpha_{i-1} \alpha_i}, S_i)$ and $\mathcal{T}_1 = \add(M_{ \alpha_i})$.  Since $\Omega(\mathcal{T}_0 ) \subset \mathcal{T}_0 $ and $\Omega(\mathcal{T}_1 ) \subset \mathcal{T}_1$, then $\Dell(S_1)\leq \pd(S_1) \leq \pd(M_{\alpha_1}) +1$, $\Dell(S_3)\leq \pd(S_3) \leq \pd(M_{\alpha_1}) +1$, and  $\Dell(S_0) = \pd(M_{\alpha_1})+1$. Hence $\Dell(A_n) \geq \lfloor \frac{n}{4} \rfloor $\\

Now, we will compute a bound for $\fin(A_n^{op})$. We consider the following two cases:

\begin{itemize}

\item $M\in \mod A_n^{op}$ such that all the vertices in $\supp(\topp(M))$ belong to the set $\{0,1 \ldots ,k+1, k+2, k+3\}$.

\item $M\in \mod A_n^{op}$ such that some of the vertices in $\supp(\topp(M))$ belong to the set $\{k+4, k+5, \ldots, n-1, n\}$.

\end{itemize}

In the first case, if $\pd (M) < \infty$, then $\pd(M) \leq \lfloor \frac{k}{4} \rfloor + 1$.\\

Since the projective $A_n^{op}$-modules $\{ P_k$, $P_{k+1}$, $P_{k+2}$, $P_{k+3}\}$ have the following vector dimensions, respectively
$$\begin{array}{cccc}
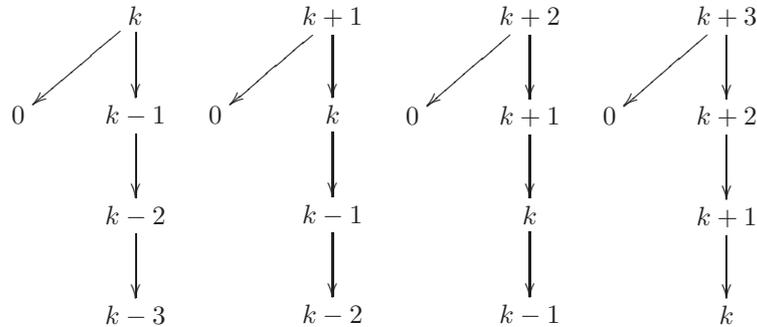
\xymatrix{ & k \ar[d] \ar[dl] \\ 0 & k-1 \ar[d] \\  & k-2  \ar[d] \\ & k-3 } & \xymatrix{ & k+1 \ar[d] \ar[dl] \\ 0 & k \ar[d] \\  & k-1  \ar[d] \\ & k-2 } & \xymatrix{ & k+2 \ar[d] \ar[dl] \\ 0 & k+1 \ar[d] \\  & k  \ar[d] \\ & k-1 } & \xymatrix{ & k+3 \ar[d] \ar[dl] \\ 0 & k+2 \ar[d] \\  & k+1  \ar[d] \\ & k } 
\end{array}$$
the $A^{op}$ modules in the second case have infinite projective dimension.
\end{ej}


\begin{thebibliography}{10}

\bibitem{ASS} I. Assem, D. Simson, A. Skowroński, {\em Elements of the representation Theory of Associative Algebras volume 1}, Cambridge University Press, Cambridge (2006).


\bibitem{BH-ZT} E. Babson, B. Huisgen-Zimmermann, R. Thomas,
{\em Generic representation theory of quivers with relations},
J. Algebra {\bf 322} (6), pp. 1877-1918 (2009).

\bibitem{BLMV} D. Bravo, M. Lanzilotta, O. Mendoza, J. Vivero. {\em Generalized Igusa-Todorov functions and Lat-Igusa-Todorov algebras}  Journal of Algebra. {\bf 580}, pp. 63-83 (2021). 

\bibitem{BM2}  M. Barrios, G. Mata, {\em On Lat-Igusa-Todorov algebras},
São Paulo J. Math. Sci. {\bf 16} (2), pp. 693-711 (2022).
 
\bibitem{BM1} M. Barrios, G. Mata, {\em On algebras of $\Omega^{n}$-finite and $\Omega^{\infty}$-infinite representation type}, to appear on J. Algebra Appl. (2023).

\bibitem{BM} M. Barrios, G. Mata,{ \em The Igusa-Todorov $\phi$-dimension on Morita Context Algebras}, Algebr. Represent. Theor. {\bf 26} pp. 3255-3269 (2023).


\bibitem{Ba} H. Bass, {\em Finitistic dimension and a homological generalization of semi-primary rings}, Trans. Amer. Math. Soc. {\bf 95} pp. 466-488 (1960).


\bibitem{G1} V. Gélinas, {\em The depth, the delooping level and the finitistic dimension},
Adv. Math. {\bf 394}, Paper No. 108052, 34 pp (2022).

\bibitem{G2} V. Gélinas,{\em The finitistic dimension of an Artin algebra with radical square zero}, (English summary) Proc. Amer. Math. Soc. {\bf 149} no.12, pp. 5001-5012 (2021).

\bibitem{GI} R. Guo, K. Igusa {\em Derived delooping levels and finitistic dimension}, arXiv:2311.00661[math.RT]

\bibitem{HLM} F. Huard, M. Lanzilotta, O. Mendoza, {\em An approach to the finitistic dimension conjecture}, J. Algebra {\bf 319} (9), pp. 3916-3934 (2008).

\bibitem{H-Z}  B. Huisgen-Zimmermann, {\em Homological domino effects and the first Finitistic Dimension Conjecture}, Invent Math  {\bf 108}, pp. 369-383 (1992).

\bibitem{IT} K. Igusa, G. Todorov, {\em On finitistic global dimension conjecture for artin algebras}, Representations of algebras and related topics, Fields Inst. Commun., {\bf 45}, American Mathematical Society, pp. 201-204 (2005).

\bibitem{LM} M. Lanzilotta, G. Mata, {\em Igusa-Todorov functions for Artin algebras}, J. Pure Appl. Algebra {\bf 222} (1), pp. 202-212 (2018). 

\bibitem{LMM} M. Lanzilotta, E. Marcos, G. Mata, {\em Igusa-Todorov functions for radical square zero algebras}, J. Algebra {\bf 487}, pp. 357-385 (2017).

\bibitem{M} H. Mochizuki, {\em Finitistic global dimension for rings}, Pac. J. Math. {\bf 15}, pp. 249-258 (1965)

\bibitem{KR} L. Kershaw, J. Rickard, {\em A finite dimensional algebra with infinite delooping level}, Ann. Repr. Th. {\bf 1} (1), pp. 61-65 (2024).

\bibitem{W} J. Wei,  {\em Finitistic dimension and Igusa-Todorov algebras}
Adv. Math. {\bf 222} (6), 2215-2226 (2009).

\end{thebibliography}
\end{document}